\documentclass{amsart}

\usepackage{amsmath,amssymb,amscd,latexsym,enumitem,amsthm, hyperref, bbding}

\usepackage[all]{xy}

\newcommand{\N}{{\mathbb{N}}}
\newcommand{\dimnuc}{\mathrm{dim}_{\mathrm{nuc}}}

\newcounter{number}[section]

\newenvironment{nummer}{\refstepcounter{number}{\bf \noindent\arabic{section}.\arabic{number}}}{}

\newcommand{\bn}{\noindent \begin{nummer} \rm}
\newcommand{\en}{\end{nummer}}

\newenvironment{ntheorem}{\noindent {\bf Theorem:} \it}{}
\newenvironment{nlemma}{\noindent {\bf Lemma:} \it}{}
\newenvironment{nprop}{\noindent {\bf Proposition:} \it}{}
\newenvironment{ndefn}{\noindent {\bf Definition:} \it}{}
\newenvironment{ncor}{\noindent {\bf Corollary:} \it}{}

\newenvironment{nquestion}{\noindent {\bf Question:} }{}
\newenvironment{nproof}{\noindent {\emph{Proof.}} }{\mbox{}\hfill 
\rule[-.2ex]{.25em}{1.8ex}}

\title[Nuclear dimension and Smale space $\mathrm{C}^*$-algebras]{Nuclear dimension and classification of $\mathrm{C}^*$-algebras associated to Smale Spaces}
\author{Robin J. Deeley}
\address{Robin J. Deeley,   Department of Mathematics, University of Hawaii,
2565 McCarthy Mall, Keller 401A, Honolulu HI 96822}
\author{Karen R. Strung}
\address{Karen R. Strung, Instytut Matematyczny Polskiej Akademii Nauk, ul. \'{S}niadeckich 8, 00-656 Warszawa, Poland}
\thanks{The first listed author was supported by ANR Project SingStar. The second listed author is supported by an IMPACT fellowship cofunded by Ministry of Science and Higher Education grant 3038/7.PR/2014/2 and EC grant PCOFUND-GA-2012-600415, and the Sonata 9 NCN grant 2015/17/D/ST1/02529}
\begin{document}

\begin{abstract}
We show that the homoclinic $\mathrm{C}^*$-algebras of mixing Smale spaces are classifiable by the Elliott invariant. To obtain this result, we prove that  the stable, unstable, and homoclinic $\mathrm{C}^*$-algebras associated to such Smale spaces have finite nuclear dimension. Our proof of finite nuclear dimension relies on Guentner, Willett, and Yu's notion of dynamic asymptotic dimension. 
\end{abstract}

\maketitle

\section{Introduction}
The classification programme for separable simple unital nuclear \mbox{$\mathrm{C}^*$-algebras} has made remarkable progress in recent years. The programme seeks to classify such \mbox{$\mathrm{C}^*$-algebras} by an invariant, known as the Elliott invariant, consisting of $K$-theory, the tracial state space, and a pairing between these objects. Since the early days of classification, there have been many interactions with the field of topological dynamical systems.  Recent research has focussed on proving the presence (or absence) of certain regularity properties of a given class of $\mathrm{C}^*$-algebras which might make it amenable to classification. In particular, a certain dimensionality condition---finiteness of the nuclear dimension---is required to avoid certain pathologies which might make a $\mathrm{C}^*$-algebra ``unclassifiable'' by Elliott invariants, that is, one may find a second $\mathrm{C}^*$-algebra with the same invariant which is nonisomorphic to the original. For a good discussion on the nuclear dimension and its relation to classification, see \cite{WinterZac:dimnuc}, \cite{SatWhiWin:NucDimZ}, or \cite{TikWhiWin:QD}. 

Minimal dynamical systems on a compact metric space $X$ (systems in which the orbit of every point under the homeomorphism $\varphi$ is dense in $X$) have provided a wealth of examples of simple separable unital nuclear $\mathrm{C}^*$-algebras via the crossed product construction. The question of their classification has seen a lot of interest, for example \cite{GioPutSkau:orbit, Put:TopStabCantor, LinPhi:MinHom, TomsWinter:PNAS, TomsWinter:minhom, Str:XxSn} to name but a few. The case of a system with mean dimension zero (see \cite{LinWeiss:MTD}) was recently fully resolved by Huaxin Lin \cite{Lin:MinDyn}. Going beyond minimal homeomorphisms, various results relating a dynamical system to an associated $\mathrm{C}^*$-algebra have been established, often focussing on crossed product $\mathrm{C}^*$-algebras \cite{KerrNow:ResFinActions, Sza:RokDimZd, RorSier:PI}, but also for other $\mathrm{C}^*$-algebras stemming from dynamics, for example, the purely infinite Cuntz--Krieger algebras associated to subshifts of finite type.

In this paper we focus on $\mathrm{C}^*$-algebras associated to certain hyperbolic dynamical systems.  A Smale space is a dynamical system $(X, \varphi)$ where $X$ is a compact metric space and $\varphi$ is a homeomorphism, with a particularly tractable local structure: at each point $x \in X$ there is a small neighbourhood which splits into two sets, which we think of as local coordinates. Along one coordinate, the systems is expanding; along the other it is contracting.  Ruelle defined Smale spaces in \cite{Rue:ThermForm} to model the restriction of an Axiom A diffeomorphism to its nonwandering set or one of its basic sets. As such this class of dynamical systems is quite diverse: it includes subshifts of finite type, William's attractors (for example solenoids), Anosov diffeomorphisms (for example hyperbolic toral automorphisms), among many others. 

The construction of $\mathrm{C}^*$-algebras associated to a Smale space is considered in \cite{Rue:NCAlgs, Put:C*Smale, PutSpi:Smale}. We briefly review these constructions. Irreducible (see Definition~\ref{irred}) Smale spaces have dense periodic points. As such, the crossed product construction, which can be seen as coming from the orbit equivalence relation (an \'etale groupoid) of a dynamical system, will be far from simple. However, when considering a Smale space, we are rather more interested in the behaviour of points at infinity (in either direction), that is, the expanding and contracting behaviour of our system. From this point of view, there are three naturally associated groupoids called the stable, unstable, and homoclinic groupoids. Suitably defined (see the preamble of Section~\ref{Smale groupoids}) these are all \'etale groupoids and using the groupoid $\mathrm{C}^*$-algebra construction as in \cite{Ren:groupoid}, result in separable, nuclear $\mathrm{C}^*$-algebras. Each of these algebras is stably finite and the homoclinic $\mathrm{C}^*$-algebra is unital. Furthermore, in the case that the original Smale space is mixing (see Definition~\ref{irred}), each of these three $\mathrm{C}^*$-algebras is simple, thus providing an excellent class of examples to examine from the point of view of the classification programme.  

In this paper we show that, given a mixing Smale space $(X, \varphi)$ the associated stable, unstable and homoclinic $\mathrm{C}^*$-algebras have finite nuclear dimension. Our proof that the Smale space $\mathrm{C}^*$-algebras have finite nuclear dimension uses the notion of dynamic asymptotic dimension for groupoids, introduced in \cite{GueWilYu:DAD}. In fact, using Smale's decomposition theorem, these results also hold under the weaker condition that the Smale space is irreducible. For the homoclinic $\mathrm{C}^*$-algebra, which is unital and has a unique tracial state, we then can deduce that in the mixing case such $\mathrm{C}^*$-algebras are distinguished by their Elliott invariants. A number of the special cases of Smale space $\mathrm{C}^*$-algebras have been studied in the context of $\mathrm{C}^*$-classification. For example, Thomsen has studied inductive limit structures \cite{Tho:HomHetExpDyn, Tho:HomHetOneSol},  in particular see \cite[Section 4.5]{Tho:HomHetExpDyn} and \cite[Introduction]{Tho:HomHetOneSol}. Tikuisis, White, and Winter have posted a preprint on quasidiagonality of nuclear $\mathrm{C}^*$-algebras \cite{TikWhiWin:QD} which implies that once we show the homoclinic algebra has finite nuclear dimension, quasidiagonality is automatic. Nevertheless, we give a proof which is independent of their result.

The stable, unstable, and homoclinic algebras each have a natural action induced from the original homeomorphism $\varphi$. While not the subject of this paper, the associated crossed product algebras, called the Ruelle algebras, are also of great interest. When the Smale space is mixing, the stable and unstable Ruelle algebras are simple and purely infinite \cite[Theorems 1.4 and 1.5]{PutSpi:Smale}. Thus, the Kirchberg--Phillips classification theorem implies they are classified by $K$-theory. When the Smale space in question is a shift of finite type, these algebras are (stably) isomorphic to the well-known Cuntz--Krieger algebras \cite{CuntzKri:OA}; in this sense the Ruelle algebras are generalisations of Cuntz--Krieger algebras. Furthermore, the Smale space $\mathrm{C}^*$-algebras formed from the stable and unstable relations play a role analogous to that of the AF-cores (studied independently by Krieger \cite{Kri:DimFun}) in the theory of Cuntz--Krieger algebras.

The paper is structured as follows. In Section~\ref{prelim} we recall some definitions about topological dynamical systems and topological groupoids and then introduce Smale spaces and their associated groupoids. In Section~\ref{dld sec} we review the notion of dynamic asymptotic dimension. We prove that the groupoids associated to an irreducible Smale space have finite dynamic asymptotic dimension in Section~\ref{dld Smale}, and hence their $\mathrm{C}^*$-algebras have finite nuclear dimension. We then prove that the homoclinic algebra is quasidiagonal in Section~\ref{QD}. Also in Section~\ref{QD}, we give further results for the homoclinic algebras of mixing Smale spaces, in particular, we show this class of $\mathrm{C}^*$-algebras is classified by the Elliott invariant. We also obtain results about the stable and unstable algebras when they contain nonzero projections.

\section{Preliminaries} \label{prelim}

We begin by recalling a few definitions about topological dynamical systems. Throughout, $X$ is an infinite compact metric space and $\varphi: X \rightarrow X$ is a homeomorphism.

\bn
\begin{ndefn}
Let $(X, \varphi)$ be a dynamical system. A point $x \in X$ is called a periodic point if there exists $n \in \mathbb{Z}_{>0}$ such that $\varphi^n(x) = x$.  A set $P$ of periodic points is $\varphi$-invariant if for any $x \in P$ we have $\varphi(x) \in P$. 
\end{ndefn}
\en

\bn
\begin{ndefn} \label{irred}
Let $(X, \varphi)$ be a dynamical system. 
\begin{enumerate}
\item $(X, \varphi)$ is called irreducible if for every ordered pair of nonempty open sets $U, V \subset X$ there exists $n \in \mathbb{Z}_{>0}$ such that $\varphi^n(U) \cap V$ is nonempty.
\item $(X, \varphi)$ is called mixing if for every ordered pair of nonempty open sets $U, V \subset X$ there exists $N \in \mathbb{Z}_{>0}$ such that $\varphi^n(U) \cap V$ is nonempty for every $n \geq N$.
\end{enumerate}
\end{ndefn}
\en

Note that a mixing dynamical system is irreducible, but the converse need not be true.

\subsection{Smale spaces}

Smale spaces are hyperbolic dynamical systems with a particularly nice local structure. They were defined by Ruelle \cite{Rue:ThermForm} based on Smale's work on Axiom A diffeomorphisms \cite{Sma:DiffDynSys}. Around every point we can find a small neighbourhood which can be decomposed into stable and unstable coordinates. The precise definition, given below, may appear abstruse, so for a good heuristic description we encourage the reader to see \cite[Section 2]{Put:C*Smale} or \cite[Section 2.1]{Put:HomThe}.

\bn
\begin{ndefn}\cite[Section 7.1]{Rue:ThermForm} \label{DefSmaleSpace}
Let $(X, d)$ be a compact metric space and let $\varphi : X \to X$ be a homeomorphism. The dynamical system $(X, \varphi)$ is called a Smale space if there are two constants $\epsilon_X > 0$ and $0<\lambda_X<1$ and a map, called the bracket map,  
\[ [ \cdot, \cdot] : X \times X \to X \]
which is defined for $x, y \in X$ such that $d(x,y)< \epsilon_X$. The bracket map is required to satisify the following axioms: 
\begin{itemize}
\item[B1.] $\left[ x, x \right] = x$, 
\item[B2.] $\left[ x, [ y, z] \right] = [ x, z]$,
\item[B3.] $\left[ [ x, y], z \right] = [ x,z ]$,
\item[B4.] $\varphi[x, y] = [ \varphi(x), \varphi(y)]$;
\end{itemize}
for $x$, $y$, and $z$ whenever both sides in the above equations are well-defined. Furthermore, the system also satisfies
\begin{itemize}
\item[C1.] For $x,y \in X$ such that $[x,y]=y$, we have $d(\varphi(x),\varphi(y)) \leq \lambda_X d(x,y)$ and
\item[C2.] For $x,y \in X$ such that $[x,y]=x$, we have $d(\varphi^{-1}(x),\varphi^{-1}(y)) \leq \lambda_X d(x,y)$.
\end{itemize}
\end{ndefn}
\en

We note that if the bracket exists, it is unique (see for example \cite[page 7]{Put:C*Smale}). 

\bn
\begin{ndefn}
Suppose $(X, \varphi)$ is a Smale space, $x \in X$, and $0<\epsilon\le \epsilon_X$ and $Y, Z \subset X$ a subset of points. Then we define the following sets
\begin{enumerate}
\item $X^s(x, \varepsilon)  :=  \left\{ y \in X \mid d(x,y) < \varepsilon, [y,x]=x \right\},$
\item $X^u(x, \varepsilon)  :=  \left\{ y \in X \mid d(x,y) < \varepsilon, [x,y]=x \right\}, $
\item $X^s(x)  :=  \left\{ y\in X \mid \lim_{n \rightarrow + \infty} d(\varphi^n(x), \varphi^n(y)) =0 \right\}, $
\item $X^u(x)   :=  \left\{ y\in X \mid \lim_{n \rightarrow - \infty} d(\varphi^n(x), \varphi^n(y)) =0 \right\},$
\item $X^s(Z) := \cup_{x\in Z} X^s(x)$ and $X^u(Z) :=\cup_{x\in Z} X^u(x)$, and
\item $X^h(Y, Z) := X^s(Y) \cap X^u(Z).$
\end{enumerate}
We say $x$ and $y$ are stably equivalent and write $x \sim_s y$ if $y \in X^s(x)$. Similarly, we say $x$ and $y$ are unstably equivalent, written $x \sim_u y$, if $y \in X^u(x)$. Points $x,y \in X$ are homoclinic if $y \in X^s(x) \cap X^u(x)$, meaning $\lim_{|n| \rightarrow \infty} d(\varphi^n(x), \varphi^n(y)) =0.$
\end{ndefn}
\en

The basic properties of $(X, \varphi)$ will not be discussed in detail here. However, we will need the following result, which is a special case of the main result of \cite{Man:ExpHomTopDim}:

\bn
\begin{nprop} \label{finCovDimSmaSpa} If $(X, \varphi)$ is a Smale space, then $X$ has finite covering dimension.
\end{nprop}
\en

 The reader can find more on the basic properties of Smale spaces in various places, for example \cite{Put:C*Smale, Put:HomThe, Rue:ThermForm}.

\subsection{Topological groupoids}

Let $G$ be a groupoid. We denote its unit space by $G^{(0)}$ and its range and source maps $r, s : G \to G^{(0)}$. The ordered pair $g, h \in G$ is composable if $s(g) = r(h)$ and their composition is denoted $gh$.  The inverse of $g \in G$ is denoted $g^{-1}$. In this paper, all groupoids are assumed to be locally compact and Hausdorff with locally compact unit space. Moreover, we will require that they are \'etale, meaning that $r$ and $s$ are local homeomorphisms.  In this case $G^{(0)}$ is an open subset of $G$ and the Haar system is given by counting measures.

To such a groupoid $G$, we associate its reduced groupoid $\mathrm{C}^*$-algebra as follows. Let $C_c(G)$ denote the (vector space of) compactly supported continuous functions on $G$. For $f_1, f_2, f \in C_c(G)$ we define multiplication and involution by
\[ (f_1 \cdot f_2) (g) = \sum_ {h_1 h_2= g} f_1(h_1) f_2(h_2), \text{ for all } g \in G \]
and 
\[ f^*(g) = \overline{f(g^{-1})}. \]
With these operations $C_c(G)$ is a $^*$-algebra. For every $x \in G^{(0)}$, let $\ell^2(s^{-1}(x))$ denote the Hilbert space of square-summable functions on $s^{-1}(x)$. From this we can define a $^*$-representation 
\[ \pi_x :  C_c(G) \to \mathcal{B}(\ell^2(s^{-1}(x))) \]
by
\[ (\pi_x(f)\xi)(g) = \sum_{h_1 h_2 = g} f(h_1) \xi(h_2), \]
for $f \in C_c(G)$, $\xi \in \ell^2(s^{-1}(x))$, $g \in s^{-1}(x)$.  The reduced groupoid $\mathrm{C}^*$-algebra, denoted $C^*_r(G)$ is the completion of $C_c(G)$ with respected to the norm
\[ \|f \| = \sup_{x \in G^{(0)}} \| \pi_x(f) \|.\]

One may also define a full groupoid $\mathrm{C}^*$-algebra $C^*(G)$. In the case that the groupoid is amenable, as is the case for the groupoids associated to a Smale space, the two coincide. 

\subsection{Groupoids associated to Smale spaces} \label{Smale groupoids} Following \cite{Rue:NCAlgs} (also see \cite{Put:C*Smale}), the homoclinic groupoid of a Smale space is defined directly using the homoclinic equivalence relation. We denote this groupoid by $G_H$; it is \'etale. 

One could also consider the groupoids given by the stable and unstable equivalence relations, but they are not \'etale. However, following the construction in \cite{PutSpi:Smale}, one can (by restricting in each case to an abstract transversal) obtain groupoids which are equivalent to the groupoids of stable and unstable equivalence (in the sense of \cite{MuhRenWil:Eqiv}) but which are in fact \'etale.  

Let $(X, \varphi)$ be an irreducible Smale space. We define the following groupoids associated to its stable and unstable equivalence relations \cite{PutSpi:Smale} (see also \cite{KPW:PDSmale, KilPut:TraceAsy, Whi:thesis}).

\bn
\begin{ndefn} \label{DefStableUnstableGroupoid} Let $P$ and $Q$ be finite $\varphi$-invariant sets of periodic points of $(X, \varphi)$. Define 
$$G_S(P):= \{ (x, y) \in X \times X  \mid x \sim_s y \hbox{ and } x, y \in X^u(P)\}.$$ 
and
$$G_U(Q):= \{ (x, y) \in X \times X  \mid x \sim_u y \hbox{ and } x, y \in X^s(Q)\}.$$ 
\end{ndefn}
\en

\bn \label{TopoFacts} To define a topology on $G_S(P)$, we follow \cite{PutSpi:Smale} (also see \cite[Section 2.2.4]{Kil:thesis}). 

Let $(x,y) \in G_S(P)$. Then, $x \sim_s y$ and hence there exists $N>0$ such that $d(\varphi^N(v), \varphi^N(w)) < \epsilon_X/2$. Now by continuity there is $0 < \delta < \epsilon_X/2$ such that, for $0 \leq n \leq N$ we have
$$ \varphi^n(X^u(w, \delta)) \subseteq X^u \left( \varphi^n(w), \epsilon_X/2 \right) $$
and 
\[ \varphi^n(X^u(v, \delta)) \subset X^u(\varphi^n(v), \epsilon_X/2). \]
From this we define a map $h: X^u(w, \delta) \rightarrow X^u(v, \delta)$ via
$$ z \mapsto \varphi^{-N}( [ \varphi^N(z), \varphi^N(v) ]). $$

By \cite{PutSpi:Smale}, $h: X^u(w, \delta) \rightarrow X^u(v, \delta)$  is a local homeomorphism, mapping $X^u(w, \delta)$ homeomorphically to a neighbourhood of $v$.

For such a 5-tuple $v, w, \delta, h, N$ as above, we define an open set 
\[ V(v,w,\delta,h,N):= \{(h(z), z) \mid z\in X^u(w, \delta)   \} \subset G_S(P).\]
 Such sets are the basic sets generating a topology for $G_S(P)$. 
\en

\bn
\begin{ntheorem}(see for example \cite[Theorem 2.17]{Kil:thesis}) \label{sumOfPropBasicSetGroupoidTop} 
We have the following properties of $G_S(P)$ and the basic sets introduced in the previous paragraph.  
\begin{enumerate}
\item The map $h$ is a local homeomorphism;
\item $V(v, w, \delta, h, N)$ gives a neighbourhood base for a topology on $G_S(P)$; 
\item $G_S(P)$ is an \'etale groupoid when we use this topology;
\item the unit space of $G_S(P)$ is $X^u(P)$; it is locally compact, but not compact.
\end{enumerate}
\end{ntheorem}
\en

The topology on $G_U(Q)$ is completely analogous and hence the details are omitted.

It is clear that $\varphi \times \varphi : X \times X \to X \times X$ defines an automorphism of the groupoid $G_H$. Moreover, since we have chosen $\varphi$-invariant sets $P$ and $Q$, it also defines an automorphism on each of $G_S(P)$ and $G_S(Q)$

\bn \label{represent}
It is shown in \cite[Theorem 1.1]{PutSpi:Smale} that each of these groupoids is amenable. In this case, the completion of any faithful $*$-representation of the compactly supported functions on the groupoid will be $^*$-isomorphic to the reduced and full groupoid $\mathrm{C}^*$-algebras. In particular, if $P$ and $Q$ are two sets of periodic points, we may represent both $C_c(G_S(Q))$ and $C_c(G_U(P))$ on the Hilbert space \[\mathcal{H} = \ell^2(X^h(P,Q))\]
where $X^h(P,Q)$ is the countable set of points in $X$ which are both stably equivalent to a point in $P$ and unstably equivalent to a point in $Q$. For further details on this construction see \cite[Sections 3 and 6]{KPW:PDSmale} or \cite[Section 1.3]{KilPut:TraceAsy}. The choice of the set $P$ (or $Q$) only affects the $\mathrm{C}^*$-algebra up to stable isomorphism.  

\begin{ndefn}
Let $P$ and $Q$ be a finite $\varphi$-invariant sets of periodic points of $(X, \varphi)$. The stable $\mathrm{C}^*$-algebra, $C^*(G_S(P))$, is defined to be the groupoid $\mathrm{C}^*$-algebra associated to $G_S(P)$ and the unstable algebra, $C^*(G_U(Q))$, is defined to be the groupoid $\mathrm{C}^*$-algebra associated to $G_U(Q)$. The homoclinic $\mathrm{C}^*$-algebra of $(X, \varphi)$ is the $\mathrm{C}^*$-algebra of the groupoid $G_H$; it is denoted by $C^*(G_H)$.
\end{ndefn}
\en

\section{Nuclear dimension and Dynamic asymptotic dimension} \label{dld sec}

We review the definitions of nuclear dimension and dynamic asymptotic dimension \cite{GueWilYu:DAD} as well as with their connection. 

\bn
Let $A$ and $B$ be $\mathrm{C}^*$-algebras. A completely positive (c.p.) map $\phi : A \to B$ is said to be {order zero} if it preserves orthogonality, that is, if $a, b \in A$ satisfy $ab = 0$ then $\phi(a)\phi(b) = 0$.

\begin{ndefn} \cite[Definition 2.1]{WinterZac:dimnuc} \label{nuc dim def} Let $A$ be a $\mathrm{C}^*$-algebra. We say that $A$ has nuclear dimension at most $n$, written $\dimnuc(A) \leq n$, if there exists a net $(F_{\lambda}, \psi_{\lambda}, \phi_{\lambda})_{\lambda \in \Lambda}$ where $F_{\lambda}$ are finite-dimensional $\mathrm{C}^*$-algebras, $\psi_{\lambda} : A \to F_{\lambda}$ and $\phi_{\lambda} : F_{\lambda} \to A$ are c.p. maps satisfying the following:
 \begin{enumerate} 
 \item $\psi_{\lambda}$ is contractive,
 \item for each $\lambda \in \Lambda$, $F_{\lambda}$ decomposes into $n+1$ ideals $F_{\lambda} = F_{\lambda}^{(0)} \oplus \cdots \oplus F_{\lambda}^{(n)}$ such that $\phi_{\lambda}|_{F_{\lambda}^{(i)}}$ is c.p.c. order zero for $i \in \{0, \dots, n\}$,
 \item $\phi_{\lambda}\circ \psi_{\lambda}(a)  \to a$ uniformly on finite subsets of $A$.
 \end{enumerate}
\end{ndefn}
If no such $n$ exists, then $A$ is said to have infinite nuclear dimension, $\dimnuc(A) = \infty$.
\en

The nuclear dimension should be thought of as a noncommutative analogue of topological covering dimension. Indeed, in the case of a commutative $\mathrm{C}^*$-algebra $\mathcal{C}_0(X)$ where $X$ is a locally compact Hausdorff space, the nuclear dimension is equal to the covering dimension of $X$ \cite[Proposition 2.4]{WinterZac:dimnuc}.  In the noncommutative case of course we no longer have such a correspondence. Nevertheless, as $\mathrm{C}^*$-algebras are often obtained by underlying structures such as dynamical systems or groupoids where one can ask if a notion of dimension, again in analogy to covering dimension, can be defined and related to the nuclear dimension of the associated $\mathrm{C}^*$-algebra. Here we deal with groupoid $\mathrm{C}^*$-algebras and the appropriate notion of dimension is the dynamic asymptotic dimension defined in \cite{GueWilYu:DAD}.

Let $G$ be an \'etale groupoid with locally compact unit space $G^{(0)}$. For any set $K \subset G$, we set 
\[ s(K) := \{ s(g) \mid g \in K\} \subset G^{(0)}, \]
and similarly
\[ r(K) := \{r(g) \mid g \in K\} \subset G^{(0)}.\]

\bn
\begin{ndefn}\cite{GueWilYu:DAD} Let $G$ be an \'etale groupoid with locally compact unit space $G^{(0)}$. We say that $G$ has dynamic asymptotic dimension less than or equal to $d$ if, for every open relatively compact $K \subset G$, there exist open sets $U_0, \dots, U_d \subset G^{(0)}$ satisfying the following:
\begin{enumerate}
\item $\{ U_0, \dots, U_d \}$ covers $s(K) \cup r(K)$ and  
\item for every $i = 0 ,\dots, d$, the groupoid generated by $\{ g \in K \mid s(g), r(g) \in U_i\}$ is a relatively compact subgroupoid of $G$.
\end{enumerate}
\end{ndefn}
\en

The groupoids associated to a Smale space are second countable and are also principal (or free, in the language of \cite{GueWilYu:DAD}). Recall that a groupoid $G$ is principal if for every $x \in G^{(0)}$ the isotropy group
\[ G^x_x := \{ g \in G \mid s(g) = r(g) = x \} \]
is trivial. The groupoid associated to an equivalence relation is always principal \cite[Example 1.2 c]{Ren:groupoid}.

The following is a special case of \cite[Theorem 8.6]{GueWilYu:DAD}.

\bn
\begin{ntheorem} \label{ndim}
Let $G$ be a second countable principal \'etale groupoid with locally compact unit space $G^{(0)}$. Furthermore assume that $G$ has dynamic asymptotic dimension at most $d$. Then 
\[ \dimnuc(\mathrm{C}_r^*(G)) \leq (d+1) (\dim G^{(0)} +1 ) -1,\]
where $\dim G^{(0)}$ is the topological covering dimension of $G^{(0)}$.
\end{ntheorem}

\en

\section{Dynamic asymptotic dimension of Smale spaces and nuclear dimension of their $\mathrm{C}^*$-algebras} \label{dld Smale}

Throughout this section we let $(X, \varphi)$ be an irreducible Smale space where the space $X$ has $\dim X = d < \infty$ and let $P$ be a finite $\varphi$-invariant set of periodic points. The stable groupoid is denoted $G_S(P)$,  as given in Definition~\ref{DefStableUnstableGroupoid}.

\bn
\begin{nlemma} \label{moveK}
Let $K \subset G_S(P)$ be an open relatively compact subset and let $U \subset G_S(P)^{(0)} = X^u(P)$ be an open subset. Suppose that there is $N \in \mathbb{N}$ such that  the subset $\{ (x,y) \in (\varphi^N \times \varphi^N)(K) \mid x, y \in U \}$ generates a relatively compact subgroupoid of $G$. Then so does $\{ (x,y) \in K \mid x, y \in \varphi^{-N}(U) \}$.
\end{nlemma}

\begin{nproof}
Since $\varphi^N \times \varphi^N$ is a homeomorphism of the groupoid, by mapping generators to generators, we see that the closures of the groupoids generated are in fact homeomorphic as subspaces of $G_S(P)$. 
\end{nproof}
\en

\bn 
\begin{nlemma} \label{pointAreCloseLemma}
Let $K \subset G_S(P)$ be an open relatively compact subset. Then for any $0< \gamma < \frac{\epsilon_X}{2}$ there exists $N_K \in \mathbb{N}$ such that, for each $n\ge N_K$, we have $x \in X^s(y, \gamma)$ whenever $(x,y) \in (\varphi^{n} \times \varphi^{n})(K)$. In particular, if $n\ge N_K$, then $d(x,y)< \gamma$ for every $(x,y) \in (\varphi^{n} \times \varphi^{n})(K)$.
\end{nlemma}
\begin{proof}
Since $\overline{K}$ is compact there exists a finite cover of $\overline{K}$ and hence $K$ by basic sets. Let $V(v_i, w_i, h_i, \delta_i, N_i)$, $i = 0, \dots, m$ be such a finite open cover. By taking the maximum over $i=0, \dots, m$, we need only prove the result for one of the basic sets. 

So let $V(v, w, h, \delta, N)$ be a basic set and $(x, y) \in V(v, w, h, \delta, N)$. By the definition of $h$, we have that $\varphi^N(x) \in X^s(\varphi^N(y), \epsilon_X)$. Using the fact that if $0 < \epsilon \le \epsilon_X$ and $x \in X^s(y, \epsilon)$, then $\varphi(x) \in X^s( \varphi(y), \lambda_X \epsilon)$ (this fact follows from B4 and C1 in Definition~\ref{DefSmaleSpace}), we have that, for any $n \geq N$, 
 \[
\varphi^n(x) \in X^s\left( \varphi^n(y), \lambda_X^{n-N} \epsilon_X \right) .
 \]
 Since $0<\lambda_X < 1$, the result follows by taking $n$ large enough so that $\lambda_X^{n-N}\epsilon_X < \gamma$. We emphasize that $n$ is independent of the choice of $(x,y) \in V(v,w,h, \delta, N)$; it only depends on $N$ and the global constants of the Smale space. 
\end{proof}
\en

\bn
\begin{nlemma} \label{NewCover}
Let $0< \gamma' \leq \gamma< \frac{\epsilon_X}{2}$ and let $K \subset G_S(P)$ be an open relatively compact subset. For any cover of $K$ by basic sets $\{ V(v_i, w_i, h_i, \delta_i, N_i)\}_{i=0}^m$, $m > 0$ with $\delta_i < \frac{\epsilon_X}{4}$,  there is $N \in \mathbb{N}$ and a cover $\mathcal{U'}$ of $K$ by basic sets $\{ V(v_i', w_i', h_i', \delta_i', N)\}_{i=0}^n$,  $n \geq m$,  satisfying
\[ 
 \varphi^N(X^u(v_i', \delta_i')) \subset X^u(\varphi^N(v_i'), \gamma) \quad\text{and}\quad  \varphi^N(X^u(w_i', \delta_i')) \subseteq X^u(\varphi^N(w_i'), \gamma)
\]
for each $0 \le i \le n$ and such that
\[ d(x,y) < \gamma' \text{ for all } (x, y) \in (\varphi^N \times \varphi^N)(K).\]
\end{nlemma}

\begin{nproof}
 By the previous lemma, for the given $\gamma'$, we may find $N_K \in \mathbb{N}$ sufficiently large so that $d(x,y) < \gamma'$ whenever $(x,y) \in (\varphi^{N_K} \times \varphi^{N_K})(K)$. Let 
 \[ N := \max \{N_0, \dots, N_m, N_K\}.\]
Fix $i$, $0\leq i \leq n$.  For each $x \in X^u(w_i, \delta_i)$, arguing as in \ref{TopoFacts}, there exists $0< \delta_x \leq \delta_i$ such that 
\[
\varphi^N(X^u(x, \delta_x)) \subseteq X^u \left( \varphi^N(x), \gamma \right).
\]
and 
\[ \varphi^N(X^u(h_i(x), \delta_x)) \subseteq X^u \left( \varphi^N(h_i(x)), \gamma \right). \]
Shrinking $\delta_x$ if necessary, we may further assume that 
\[  X^u(h_i(x), \delta_x)  \subset X^u(v_i, \delta_i). \]
Let $h_x$ denote the restriction $h_x := h_i|_{X^u(x, \delta_x)}$. Then $V(h_i(x), x, h_x, \delta_x, N)$ is a basic set. 

Now let $(h_x(z), z) \in V(h_x(x), x, \delta_x, h_x N)$. Then 
\begin{eqnarray*}
h_x(z) &=&  \varphi^{-N} ( [ \varphi^N(z),  \varphi^N(h_i(x))])\\  
&=& \varphi^{-N_i - (N - N_i))} ( [ \varphi^{-N_i - (N - N_i)}(z),  \varphi^{-N_i - (N - N_i)}(h_i(x))])\\
&=& \varphi^{-N_i - (N - N_i))} ( \varphi^{N- N_i}[ \varphi^{N_i}(z), \varphi^{N_i)}(h_i(x))])\\
&=&  \varphi^{-N_i} ( [ \varphi^{N_i}(z),  \varphi^{N_i}(h_i(x))]),
\end{eqnarray*}
by the fact that the bracket in the above computations is everywhere defined, allowing us to apply B4 of Definition~\ref{DefSmaleSpace}.

Furthermore,
\begin{eqnarray*}
[ \varphi^{N_i}(z),  \varphi^{N_i}(h_i(x))] &=& [ \varphi^{N_i}(z), [\varphi^{N_i}(h_i(x)), \varphi^{N_i}(v)] ] \\
&=& [ \varphi^{N_i}(z),  \varphi^{N_i}(v)] ,
\end{eqnarray*}
so $h_x(z) = \varphi^{-N_i}([ \varphi^{N_i}(z),  \varphi^{N_i}(v)]) = h_i(z)$ and we have that 
$$V(h_i(x), x, h_x, \delta_x, N) \subseteq V(w_i, v_i, \delta_i, h_i, N_i).$$ 
Since we may do this for each $i$ and $x$, the argument is completed using the compactness of $\overline{K}$. 
\end{nproof}
\en

For $0<\epsilon< \frac{\epsilon_X}{2}$, we let $U(w, \epsilon)$ denote the image of $X^u(w,\epsilon) \times X^s(w, \epsilon)$ under the bracket map.  

\bn
\begin{nlemma} \label{precomp}
Let $0< \gamma' \leq \gamma < \frac{\epsilon_X}{4}$ and let $K \subset G_S(P)$ be an open relatively compact subset. Then there is $N \in \mathbb{N}$ such that for every $w \in X$, the set
\[ \{(x,y) \in (\varphi^N \times \varphi^N)(K) \mid x, y \in U(w, \gamma)\} \]
generates an open relatively compact groupoid and such that
\[ d(x, y) < \gamma' \text{ for every } (x, y) \in (\varphi^N \times \varphi^N)(K). \]
\end{nlemma}

\begin{nproof}
Since $w$ and $\gamma$ are fixed throughout the proof, we set 
\[ U := U(w, \gamma).\] 
By Lemma~\ref{NewCover}, there is $N \in \mathbb{N}$ and a finite subcover of $K$ by $\{ V(v_i, w_i, h_i, \delta_i, N)\}_{i = 0}^m$ such that 
\[ 
 \varphi^N(X^u(v_i, \delta_i)) \subseteq X^u(\varphi^N(v_i), \gamma)) \quad\text{and}\quad  \varphi^N(X^u(w_i, \delta_i)) \subseteq X^u(\varphi^N(w_i), \gamma))
\]
for each $0 \le i \le m$ and such that $d(x, y) < \gamma'$ for every $(x, y) \in \varphi^N \times \varphi^N(K)$. It also follows from Lemma \ref{NewCover} that for each $0 \leq i \le m$ there is some $\delta_i' < \frac{\epsilon_X}{4}$ such that
\[ (\varphi^N \times \varphi^N)(V(v_i, w_i, h_i, \delta_i, N) ) \subseteq V(\varphi^N(v_i), \varphi^N(w_i), [ \, \cdot \,, \varphi^N(v_i)], \delta_i', 0),\]
In particular, $\{ V(\varphi^N(v_i), \varphi^N(w_i), [ \, \cdot \,, \varphi^N(v_i)], \delta_i', 0)\}_{i = 0}^m$ is a cover of $(\varphi^N \times \varphi^N)(K)$.

We need only consider $i$ such that 
\begin{equation} \label{nontrivialIntersection}
X^u(\varphi^N(w_i), \delta_i') \cap U \hbox{ and }X^u(\varphi^N(v_i), \delta_i') \cap U
\end{equation}
are both nonempty. Suppose that there are $n \leq m$ such sets. By reordering, let  $\{ V(\varphi^N(v_i),\varphi^N(w_i), [ \, \cdot \,, \varphi^N(v_i)] , \delta_i', 0)\}_{i=0}^n$ be the sets in the cover of $(\varphi^N \times \varphi^N)(K)$ with this property. 

By \eqref{nontrivialIntersection}, the triangle inequality and the choice of $\delta_i$ and $\gamma$, we may define the following points $r_0, \dots r_{2n_1} \in X$ using the bracket: 
\[ r_l = \left\{ \begin{array}{ccc}  \lbrack w, \varphi^N(v_i) \rbrack & \, & \text{for } 0 \leq l \leq n,\\
 \lbrack  w, \varphi^N(w_i)\rbrack & \, &\text{for } n+1 \leq l \leq  2n+1.\\
   \end{array} \right. \]
   Let
\[ H = \cup_{1 \leq l, m \leq 2n} V(r_l, r_m, \lbrack \, \cdot \, , r_l \rbrack, \gamma).\]

We first claim that $H$ is an open relatively compact subgroupoid $G_S(P)$. Openness and relative compactness follow because the set is the union of finitely many open relatively compact subsets of $G_S(P)$. Thus we need only show that $H$ is closed under inverses and, where defined, products. So let $(x, y) \in H$. Then $(x, y) \in V(r_l, r_m, \lbrack \, \cdot \, , r_l \rbrack, \gamma)$ for some $0 \leq l, m \leq 2n+1$, so we have that
\[ x \in X^u(r_m, \gamma) \text{ and } y \in X^u(r_l, \gamma) .
\]
and that
\[ \lbrack y , r_m \rbrack = x.\]
But also  
\[ \lbrack x , r_l \rbrack = y,\]
so we have that 
\[ (y, x) = ( \lbrack x , r_l \rbrack , x)  \in V(r_m, r_l, \lbrack \, \cdot \, , r_m\rbrack, \gamma) \subset H,\]
whence we see that $H$ is closed under inverses.
Now suppose that $(x, y), (y, z) \in H$ and that they are composable as elements of $G_S(P)$.  Then, as above,  for some $k, l, m$, we have that 
\[ x \in  X^u(r_k, \gamma), \quad y \in X^u(r_l, \gamma), \quad z \in X^u(r_m,\gamma),\] 
and moreover that 
\[ \lbrack y , r_k \rbrack = x , \text{ and } \lbrack z , r_l \rbrack = y. \]
Since $\gamma$ is sufficiently small, we may compose the bracket to get that
\[ x =  \lbrack \lbrack z , r_l \rbrack , r_k \rbrack = \lbrack z , r_k\rbrack,\]
thus
\[ (x, z) = (\lbrack z , r_k\rbrack, z) \in V_{r_m, r_k, \lbrack \, \cdot \, r_k \rbrack, \gamma} \subset H.\]

Next we claim the groupoid generated by $\{(x,y) \in (\varphi^N \times \varphi^N)(K) \mid x,y \in U \}$ is contained in $H$. To prove this claim, it is enough to show  
\[ ((\varphi^N \times \varphi^N)(K)) \cap (U \times U) \subset H.\]

Let $(x, y) \in ((\varphi^N \times \varphi^N)(K)) \cap (U\times U)$. Then  
\[ (x, y) \in V(\varphi^N(v_i), \varphi^N(w_i), [ \, \cdot \,, \varphi^N(v_i)], \delta_i')\]
for some $0 \leq i \leq n'$. Thus $x = \lbrack y, \varphi^N(v_i) \rbrack$. Since $y \in X^u(w, \gamma)$, we have that $\lbrack y, w \rbrack$ is well-defined and so too is the composition $ \lbrack  y,  \lbrack w , \varphi^N(v_i) \rbrack \rbrack$.
Thus
\[ \lbrack y, r_i \rbrack =  \lbrack  y, \lbrack w, \varphi^N(v_i) \rbrack \rbrack = \lbrack y, \varphi^N(v_i) \rbrack = x,\]
hence
\[ (x, y) = (\lbrack y, r_i \rbrack, y) \in H .\] 
It now follows from the two claims that $\{(x,y) \in (\varphi^N \times \varphi^N)(K) \mid x,y \in U \}$ generates a relatively compact subgroupoid of $G_S(P)$.
\end{nproof}
\en

\bn
\begin{nlemma} \label{spaced cover}
Let $(X, \varphi)$ be a Smale space (where $X$ is an infinite compact metric space with $\dim X = d < \infty$) and $0< \gamma < \frac{\epsilon_X}{4}$. Then, there exists finite open covers of $X$, $\{ U_j \}_{j=1}^M$ and $\{ V_j \}_{j=1}^M$, with the following properties:
\begin{enumerate}
\item for each $j$, there exists $w\in X$ such that $U_j \subseteq U(w, \gamma)$;
\item the cover $\{ U_j\}_{j=1}^M$ has order $d$;
\item for each $j$, $V_j \subseteq U_j$;
\item for each $j_1$, $j_2$, $U_{j_1} \cap U_{j_2} = \emptyset$ implies $\overline{V}_{j_1} \cap \overline{V}_{j_2} = \emptyset$;
\item there exists $\gamma' >0$, such that for any $j_1$, $j_2$, $x\in V_{j_1}$, $y\in V_{j_2}$ with $U_{j_1} \cap U_{j_2}=\emptyset$, we have that $d(x,y)> \gamma'$. Furthermore, we can assume $\gamma' < \gamma$.
\end{enumerate}
\end{nlemma}
\begin{proof}
By compactness, there is a finite open cover of $X$ by sets of the form $U(w, \gamma)$, $w \in X$. The definition of covering dimension applied to $X$ implies that there is a finite refinement of that cover to a cover of order $d$. We denote the refined cover by $\{U_i\}_{i=1}^M$. We construct the second cover inductively on $i=1, \ldots, M$. 

For $i=1$, we take $V_1=U_1$.

Next suppose $\{V_1, \ldots, V_i \}$ is a set of open subset of $X$ such that
\begin{enumerate}
\item $\{V_1, \ldots, V_i \} \cup \{ U_j \}_{j=i+1}^M$ is an open cover of $X$;
\item $V_s \subseteq U_s$ for each $1\le s \le i$;
\item for each $j_1$ and $j_2$ in $\{1, \ldots, i\}$, whenever $U_{j_1} \cap U_{j_2} = \emptyset$, $\overline{V}_{j_1}\cap \overline{V}_{j_2}=\emptyset$.
\end{enumerate}
Let $J_{i+1} \subseteq \{ 1, \ldots, M \}$ be the set of indices $k$ such that $U_{i+1} \cap U_k= \emptyset$, but $\overline{U}_{i+1} \cap \overline{U}_k \neq \emptyset$. Then $F_{i+1}= \cup_{k \in J_{i+1}} \overline{U}_{i+1} \cap \overline{U}_k$ is closed and hence compact. Furthermore, since $\{V_1, \ldots, V_i \} \cup \{ U_j \}_{j=i+1}^M$ covers $X$ and $F_{i+1} \cap U_{i+1} = \emptyset$, we have that 
\[
F_{i+1} \subseteq (\cup_{j=1}^i V_j) \cup (\cup_{j=i+2}^M U_j).
\]
Take an open set $G_{i+1}$ such that
\[
F_{i+1} \subseteq G_{i+1} \subseteq \overline{G}_{i+1} \subseteq (\cup_{j=1}^i V_j) \cup (\cup_{j=i+2}^M U_j).
\]
Let $V_{i+1}= U_{i+1} \cap \overline{G}_{i+1}^c$, which is open and contained in $U_{i+1}$. One then checks that
\begin{enumerate}
\item $\{V_1, \ldots, V_{i+1} \} \cup \{ U_j \}_{j=i+2}^M$ is an open cover of $X$;
\item $V_s \subseteq U_s$ for each $1\le s \le i+1$;
\item for each $j_1$ and $j_2$ in $\{1, \ldots, i+1\}$, whenever $U_{j_1} \cap U_{j_2} = \emptyset$, $\overline{V}_{j_1}\cap \overline{V}_{j_2}=\emptyset$.
\end{enumerate}
Thus, using induction, we have covers $\{U_i \}_{j=1}^M$ and $\{V_j \}_{j=1}^M$ that satisfy items (i) to (iv) in the statement of the lemma. 

We show that we also have item (v). For $j, j' \in \{ 1, \ldots M \}$ with $U_j \cap U_{j'}= \emptyset$, let $\gamma_{j, j'}$ be a positive constant such that $d(x,y)> \gamma_{j, j'}$ whenever $x \in V_j$ and $y\in V_{j'}$. We note that $\gamma_{j, j'}$ can be taken to be positive because $\overline{V}_j \cap \overline{V}_{j'}= \emptyset$. Let 
\[
\gamma' = {\rm min}\{ \gamma_{j, j'} \mid j, j' \hbox{ such that }U_j \cap U_{j'}=\emptyset \}
\]
By construction, $\gamma'$ (and any positive constant less than or equal to it) has the required property.
\end{proof}
\en

\bn
\begin{nlemma} \label{precompForVs}
Suppose $0< \gamma < \frac{\epsilon_X}{4}$ and $\{ U_j \}_{j=1}^M$, $\{V_j\}_{j=1}^M$, and $\gamma'$ are as in the statement of Lemma \ref{spaced cover}. Furthermore, suppose $K$ is an open relatively compact subset of $G_S(P)$ and $N\in \N$ is as in the statement of Lemma \ref{precomp} (it depends on $0< \gamma' < \gamma< \frac{\epsilon_X}{4}$ along with $K$). Then, for each $j$, the groupoid generated by 
\[
\{(x,y) \in (\varphi^N \times \varphi^N)(K) \mid x, y \in V_j \}
\]
is relatively compact.
\end{nlemma}
\begin{proof}
By assumption, for each $j$, there exists $w\in X$ such that $V_j \subseteq U_j \subseteq U(w, \gamma)$. The required result then follows from Lemma \ref{precomp}.
\end{proof}
\en

\bn
\begin{ntheorem} \label{dld}
Let $(X, \varphi)$ be a irreducible Smale space and $P$ a finite $\varphi$-invariant set of points. Then both the stable groupoid $G_S(P)$ and the unstable groupoid $G_U(P)$ have dynamic asymptotic dimension at most $\dim X$.
\end{ntheorem}

\begin{nproof}
By Proposition \ref{finCovDimSmaSpa}, $X$ has finite covering dimension; set $d = \dim X$. Fix $0< \gamma < \frac{\epsilon_X}{4}$ and take $\{ U_j \}_{j=1}^M$, $\{V_j\}_{j=1}^M$, and $\gamma'$ as in Lemma \ref{spaced cover}. Furthermore, let $K \subseteq G_S(P)$ be an open relatively compact subset.

Take $N\in\mathbb{N}$ as in the statement of Lemma~\ref{precomp} and let $\tilde{K}= (\varphi^N \times \varphi^N)(K)$. Note that, in particular, $d(x, y) < \gamma'$ whenever $x,y \in \tilde{K}$. It follows from  Lemma~\ref{moveK} that it is enough to show that there are open sets $\{\tilde{V}_0, \dots, \tilde{V}_d\}$ which form an open cover of $s(\tilde{K}) \cup r(\tilde{K})$ and such that, for each $0 \leq i \leq d$, the set
\[ \{ (x,y) \in \tilde{K} \mid x, y \in \tilde{V}_i \} \]
generates a relatively compact subgroupoid.

Let $V$ be an open relatively compact subset of $X^u(P)$ that contains $s(\tilde{K})\cup r(\tilde{K})$. Then, for each $1\le j \le M$, $U_j \cap V$ and $V_j\cap V$ are open subsets of $X^u(P)$. Moreover, both $\{ U_j\cap V \}_{j=1}^M$ and $\{V_j \cap V\}_{j=1}^M$ are covers of $s(\tilde{K}) \cup r(\tilde{K})$.

Since $\{ U_j \}_{j=1}^M$ has order $d$, we may decompose $\{ 1, 2, \ldots, M \}$ into subsets $J_0, \dots, J_d$ such that, for every $i = 0, \dots, d$, the sets $U_j$, $j \in J_i$ are pairwise disjoint. 

For every $i = 0, \dots, d$, let
\[ \tilde{V}_i = \cup_{j \in J_i} V_j \cap V. \]
Then $\{\tilde{V}_0, \dots, \tilde{V}_d\}$ is an open cover of $s(\tilde{K}) \cup r(\tilde{K})$. 

Using item (v) in the statement of Lemma \ref{spaced cover} and the fact that the $U_j, j \in J_i$ are pairwise disjoint, we have that for $j \neq j' \in J_i$, $d(x, y) > \gamma'$  whenever $x \in V_j$ and $y \in V_{j'}$. It follows that, for each $i=0, \dots , d$, the groupoid generated by 
\[ \{ (x,y) \in \tilde{K} \mid x, y \in \tilde{V}_i \} \]
is simply the disjoint union (over $j\in J_i$) of the groupoids generated by 
\[ \{ (x,y) \in \tilde{K} \mid x, y \in V_j  \}. \]
Thus, using Lemma~\ref{precompForVs}, we have that for each $i = 0, \dots, d$, the set $\{ (x,y) \in (\varphi^N \times \varphi^N)(K) \mid x, y \in \tilde{V}_i \}$ generates a relatively compact groupoid. In other words, $\tilde{V}_0, \ldots, \tilde{V}_d$ is the required cover of $s(\tilde{K})\cup r(\tilde{K})$. This completes the proof for $G_S(P)$.
To show that $G_U(P)$ has dynamic asymptotic dimension at most $d$, we simply observe that $G_U(P)$ is the stable groupoid of the irreducible Smale space $(X, \varphi^{-1})$ and the result follows. 
 \end{nproof}
\en

\bn
\begin{ncor} \label{finNucDimSmaCstarThm}
Let $(X, \varphi)$ be an irreducible Smale space and $P$ a finite set of $\varphi$-invariant periodic points. Then the stable, unstable and homoclinic $\mathrm{C}^*$-algebras each have finite nuclear dimension.
\end{ncor}

\begin{nproof} 
Combining Theorem~\ref{dld} with Theorem~\ref{ndim}, we have that $\dimnuc(\mathrm{C}^*(G_S(P))$ and $\mathrm{C}^*(G_U(P))$ are bounded by $(\dim(X) -1)(\dim(X^u(P))-1)+1$ and $(\dim(X)-1)(\dim(X^s(P)-1)+1$ respectively. The result for $\mathrm{C}^*(G_H)$ follows from the fact that it is Morita equivalent to $\mathrm{C}^*(G_S(P)) \otimes \mathrm{C}^*(G_U(P))$ \cite[Theorem 3.1]{Put:C*Smale}.
\end{nproof}
\en

\section{Structure and Classification} \label{QD}
To begin, we will prove the following result:

\bn
\begin{ntheorem} \label{quasid}
Let $(X, \varphi)$ be a mixing Smale space and $P$ be a set of periodic points as in Definition \ref{DefStableUnstableGroupoid}. Then the homoclinic $\mathrm{C}^*$-algebra is quasidiagonal.
\end{ntheorem}
\en

To prove the theorem, we require the following lemmas. The first lemma is likely well-known and is implicit in \cite{BlaKir:GenIndLim}, but the proof is short so we include it.

\bn
\begin{nlemma}
Suppose $A$ is a separable, nuclear $\mathrm{C}^*$-algebra and that
$$A \subset \frac{\Pi_{n \in \N}  \mathcal{K}(\mathcal{H})}{ \bigoplus_{n \in \N} \mathcal{K}(\mathcal{H})}.$$
Then $A$ is quasidiagonal.
\end{nlemma}

\begin{nproof}
The compact operators are MF (see \cite[Definition 3.2.1]{BlaKir:GenIndLim}), so by \cite[Corollary 3.4.3]{BlaKir:GenIndLim} $A$ is also MF. Since $A$ is nuclear, this is equivalent to $A$ being quasidiagonal by \cite[Theorem 5.2.2]{BlaKir:GenIndLim}.
\end{nproof}
\en

In the proof of the next lemma, we use basic properties of the representations of $C^*(G_S(P))$ and $C^*(G_U(Q))$ on the Hilbert space $\mathcal{H} = \ell^2(X^h(P, Q))$ (see \ref{represent}). These results can be found in \cite[Section 6]{KPW:PDSmale}, \cite[Theorem 1.3]{KilPut:TraceAsy} and \cite[Section 4]{Whi:thesis}. The last of these references contains complete proofs of each of the results we require. Note that we do not make use of \cite[Lemma 4.4.9]{Whi:thesis} so we do not need to assume that the sets of periodic points used to defined the stable and unstable algebras are disjoint.

\bn
\begin{nlemma}
Let $P$ and $Q$ be finite sets of $\varphi$-invariant periodic points. Suppose that $C^*(G_S(P))$ and $C^*(G_U(Q))$ are the stable and unstable $\mathrm{C}^*$-algebras associated to a mixing Smale space. Then there exists an injective $*$-homomorphism
\[\Phi:  C^*(G_S(P)) \otimes C^*(G_U(Q)) \to \frac{\Pi_{n \in \N}  \mathcal{K}(\mathcal{H})}{ \bigoplus_{n \in \N} \mathcal{K}(\mathcal{H})}. \]
\end{nlemma}
\begin{nproof}
Let $\alpha_s : C^*(G_S(P)) \to C^*(G_S(P))$ the $^*$-automorphism induced by $\varphi$. By \cite[Lemma 6.1]{KPW:PDSmale}, for any $a \in C^*(G_S(P))$ and $b\in C^*(G_U(Q))$, $ab\in \mathcal{K}(\mathcal{H})$. Thus, the map on elementary tensors defined by $a\otimes b \mapsto (\alpha^n_s(a) b)_{n\in \N}$ extends to a linear map 
\[
\Phi : C^*(G_S(P))\otimes C^*(G_U(Q)) \rightarrow \frac{\Pi_{n \in \N}  \mathcal{K}(\mathcal{H})}{ \bigoplus_{n \in \N} \mathcal{K}(\mathcal{H})}.
\]
We show that $\Phi$ is a $^*$-homomorphism. To do so, we must show that, for each $a, a' \in C^*(G_S(P))$ and $b, b' \in C^*(G_U(Q))$, 
\begin{align}
\left( \alpha^n_s(a) b \alpha^n_s(a')b' - \alpha^n_s(a)\alpha^n_s(a')bb' \right)_{n\in \N} & \in \bigoplus_{n \in \N} \mathcal{K}(\mathcal{H}) \hbox{ and } \\
\left( \alpha^n_s(a^*)b^* - ( \alpha^n(a)b)^* \right)_{n\in \N} & \in \bigoplus_{n \in \N} \mathcal{K}(\mathcal{H}).
\end{align}
We prove only the former as the proofs are rather similar. Let $\epsilon >0$. We must show that there exists $N \in \N$ such that for all $n \ge N$
\[
\| \alpha^n_s(a) b \alpha^n_s(a')b' - \alpha^n_s(a)\alpha^n_s(a')bb' \| < \epsilon.
\]
By \cite[Lemma 6.3]{KPW:PDSmale}, there exists $N \in \N$ such that, for any $n \ge N$,
\[
\| b \alpha^n_s(a')-\alpha^n_s(a')b \| < \frac{\epsilon}{\| a \| \| b'\|}.
\]
Hence
\[
\| \alpha^n_s(a) b \alpha^n_s(a')b' - \alpha^n_s(a)\alpha^n_s(a')bb' \| \le \| a \| \|b \alpha^n_s(a')-\alpha^n_s(a')b \| \|b' \| < \epsilon.
\]
Finally, we must show that $\Phi$ is injective. Recall that $C^*(G_S(P))$ and $C^*(G_U(Q))$ are both simple so $C^*(G_S(P)) \otimes C^*(G_U(Q))$ is also simple. Thus, we need only show that $\Phi$ is nonzero on a single element in $C^*(G_S(P)) \otimes C^*(G_U(Q))$. The existence of such an element follows from the proof of Lemma 4.4.13 in \cite{Whi:thesis}. In that proof, elements $a\in C^*(G_S(P))$ and $b\in C^*(G_U(Q))$ are constructed with the property that $\| \alpha_s^{n_k}(a) b \| \geq 1$ for a strictly increasing sequence $n_k$.
\end{nproof}
\en

We can now prove the main result of the section. 

{\emph{Proof of Theorem~\ref{quasid}.\/}}  The previous two lemmas imply that the tensor product $C^*(G_S(P)) \otimes C^*(G_U(Q))$ is quasidiagonal. Since $C^*(G_H)$ is stably isomorphic to $C^*(G_S(P)) \otimes C^*(G_U(Q))$ it is also quasidiagonal \cite[Theorem 20]{Had:QD}. \qed

\bn
\begin{ncor}
Let $(X, \varphi)$ be an irreducible Smale space. Then its homoclinic $\mathrm{C}^*$-algebra is quasidiagonal.
\end{ncor}

\begin{nproof}
It follows from Smale's decomposition theorem \cite[Theorem 6.2]{Sma:DiffDynSys} that the Smale space $\mathrm{C}^*$-algebras of an irreducible Smale space have a special form (see for example \cite[Section 2.5]{Kil:thesis}): they are each the direct sum of finitely many $\mathrm{C}^*$-algebras each of which is the $\mathrm{C}^*$-algebra of a mixing Smale space. The result then follows from Theorem~\ref{quasid} because quasidiagonality is preserved under direct sums.
\end{nproof}
\en

If, in the definition of nuclear dimension (Definition \ref{nuc dim def}), we in addition require the maps $\phi_{\lambda}$ to be contractive, then we have the definition for the decomposition rank, a precursor to nuclear dimension which first appeared in \cite{KirWinter:dr}. For the homoclinic algebra we obtain the following bound on both the decomposition rank as well as a smaller bound on the nuclear dimension.

\bn\begin{ncor} \label{decRanMatSat}
Let $(X, \varphi)$ be a mixing Smale space. Then the associated homoclinic $\mathrm{C}^*$-algebra $A$ has decomposition rank at most one.
\end{ncor}

\begin{nproof}
By \cite[Theorem 1.1]{MatSat:UHF}, since $A$ is simple, separable, unital, nuclear and has unique tracial state, it has finite decomposition rank. By \cite[Corollary 8.6]{BBSTWW:2Col} this is equivalent to having decomposition rank at most one. Hence $A$ also has nuclear dimension at most one. 
\end{nproof}
\en

\bn
\begin{ncor} \label{ratTAF}
Let $A$ be the homoclinic $\mathrm{C}^*$-algebra of a mixing Smale space $(X, \varphi)$. Then $A \otimes \mathcal{U}$ is tracially approximately finite (TAF) for any UHF algebra $\mathcal{U}$ of infinite type .
\end{ncor}
\en

\begin{nproof}
The $\mathrm{C}^*$-algebra $A$ is separable, unital and nuclear and has a unique trace.  Since $(X, \varphi)$ is mixing, $A$ is simple.  By Theorem~\ref{quasid} it is also quasidiagonal. Thus it follows from \cite[Theorem 6.1]{MatSat:UHF} that $A \otimes \mathcal{U}$ is TAF for any UHF algebra $\mathcal{U}$ of infinite type.
\end{nproof}

\bn
\begin{ntheorem} \label{mainClassThm}
The homoclinic algebras associated to mixing Smale spaces are contained in a class of $\mathrm{C}^*$-algebras that is classified by the Elliott invariant. In particular,
\begin{enumerate}
\item the homoclinic algebra associated to a mixing Smale space is approximately subhomogeneous,
\item if $A$ and $B$ are the homoclinic algebras associated to the mixing Smale spaces $(X, \varphi)$ and  $(Y, \psi)$, then an isomorphism 
\[ \phi: (K_0(A), K_0(A)_+, [1_A], K_1(A)) \to (K_0(B), K_0(B)_+, [1_B], K_1(B))\]
lifts to a $^*$-isomorphism 
\[ \Phi: A \to B \]
inducing $\phi$. 
\end{enumerate}
\end{ntheorem}

\begin{nproof}
Once we know finite decomposition rank, this in fact follows directly from the main results in  \cite{EllGonLinNiu:ClaFinDecRan}, since the homoclinic groupoid of $(X, \varphi)$ is amenable and hence its $\mathrm{C}^*$-algebra satisfies the UCT \cite[Proposition 10.7]{Tu:Groupoids}. However, we do not need the full strength of this classification result and can appeal to earlier work. By Corollary \ref{ratTAF}, for any homoclinic algebra $A$ corresponding to a mixing Smale space $(X, \phi)$ we have that $A \otimes \mathcal{U}$ is tracially approximately finite for any UHF algebra $\mathcal{U}$ of infinite type.    By applying  \cite[Theorem 5.4]{LinNiu:KKlifting} the $\mathrm{C}^*$-algebra $A \otimes \mathcal{Z}$ belongs to a class classified by the Elliott invariant. But $A$ also has finite nuclear dimension by Theorem~\ref{finNucDimSmaCstarThm}, so since it is also separable, simple, unital and nonelementary, $A$ tensorially absorbs $\mathcal{Z}$ \cite[Corollary 7.3]{Win:Z-stabNucDim}. The result follows.
\end{nproof}
\en

In the case that the stable and unstable algebras associated to a mixing Smale space both have nonzero projections, we can also say something about their structure.

\bn
\begin{ntheorem} \label{stableAndUnstableWithProj}
Suppose that $(X, \varphi)$ is a mixing Smale space and that the associated unstable and stable algebras each contain a nontrivial projection. Then they are approximately subhomogeneous (ASH) algebras. 
\end{ntheorem}

\begin{nproof}
Recall that the choice of the particular sets of periodic periods $P$ and $Q$ in the definition of the stable and unstable $\mathrm{C}^*$-algebras only affects the algebras up to stable isomorphism. We prove that the stable algebra $C^*(G_S(P))$ is ASH since the proof for $C^*(G_U(Q))$ is identical upon changing the roles of the two algebras. 

Since $C^*(G_U(Q))$  contains a nonzero projection it is isomorphic to a separable subalgebra of $C^*(G_S(P))\otimes C^*(G_U(Q))$. By (the proof of) Theorem~\ref{quasid}, we have that $C^*(G_S(P))\otimes C^*(G_U(Q))$ is quasidiagonal so it follows that $C^*(G_S(P))$ is also quasidiagonal. Thus the stable algebra is quasidiagonal, satisfies the UCT (again by \cite[Proposition 10.7]{Tu:Groupoids}) and, by Corollary~\ref{finNucDimSmaCstarThm}, has finite nuclear dimension. 

Let $p$ be a nonzero projection in $C^*(G_S(P))$. Since $C^*(G_S(P))$ is simple, $p$ is full so $C^*(G_S(P))$ is stably isomorphic to $p(C^*(G_S(P)))p$ \cite[Corollary 2.6]{Bro:stabher}. Thus $C^*(G_S(P))$ is stably isomorphic to a $\mathrm{C}^*$-algebra of rational generalised tracial rank at most one by \cite[Theorem 4.3]{EllGonLinNiu:ClaFinDecRan} and \cite[Theorem A]{TikWhiWin:QD}, and hence has an ASH model given by \cite{Ell:invariant}. Since any stabilisation and any hereditary subalgebra of an ASH algebra is again ASH, it follows that $C^*(G_S(P))$ is ASH.
\end{nproof}
\en

In \cite{Wie:SmaleLimits} Wieler constructs many examples of Smale spaces with totally disconnected stable sets. It is not difficult to show that in this case the unstable algebra will have many projections; the unit space of the unstable groupoid will be totally disconnected in this case. It is natural to ask the following:

\bn
\begin{nquestion} \label{projQues}
When do the stable and unstable algebras of a mixing Smale space contain projections?
\end{nquestion}
\en

Based on examples, it is plausible that the answer to this question is ``always".

\section{Further Remarks}
Though our classification results only apply to the homoclinic $\mathrm{C}^*$-algebra, there are situations where our results for the stable and unstable algebras also have consequences for classification. In particular, the proof of Theorem \ref{mainClassThm} implies that the class of unital $\mathrm{C}^*$ algebras Morita equivalent to the unstable algebra of a Smale space is classifiable by the Elliott invariant. 

An explicit example of this situation comes from certain substitution tilings.  Let $A$ be the $\mathrm{C}^*$-algebra considered in \cite[Section 5]{KelPut:TilCstarKth}. Then $A$ is unital and Morita equivalent to the Smale space constructed by Anderson and Putnam in \cite[Section 3]{AndPut:TopInvSubTilCstar}. We should note that for this specific case, classification can also be proved using \cite[Corollary 3.2]{Win:ClassCrossProd} and the main result of \cite{SadWil:TilSpaFibBun} (also see \cite[Remark 7.10]{Sta:FinActTilAlg}). Phillips has asked if these substitution tiling $\mathrm{C}^*$-algebras are TAF (again see \cite[Remark 7.10]{Sta:FinActTilAlg}). Our results do not imply this, but classification could be used to reduce this to a question about $K$-theory. For example, if the $K$-theory is torsion-free then the relevant $\mathrm{C}^*$-algebra is A$\mathbb{T}$.

In general, our results indicate the importance of computing the $K$-theory of Smale space $\mathrm{C}^*$-algebras. This is already an active area of research. Some calculations have been done in particular cases, see for example \cite{Nekra:SelfSim, Put:C*Smale, Tho:HomHetExpDyn, Tho:HomHetOneSol}. Moreover, Putnam has constructed a homology theory for Smale spaces \cite{Put:HomThe}. A link between Putnam's homology and the $K$-theory of the $\mathrm{C}^*$-algebra has been conjectured in \cite[Section 8.4]{Put:HomThe}. Our result are relevant for the computation of the $K$-theory. This is because finite dynamic asymptotic dimension (see \cite[Introduction]{GueWilYu:DAD} for details) should lead to an in principal commutation of the $K$-theory. It will be interesting to see if this leads to an effective way to compute the $K$-theory (at least for special classes of Smale spaces).  

Finally, we observe that in certain special cases, an inductive limit structure for the homoclinic $\mathrm{C}^*$-algebra can be found directly from the dynamics. For example, in the case of a shift of finite type, the homoclinic (and also the stable and unstable) $\mathrm{C}^*$-algebra is AF. One dimensional solenoids are another example, see the work of Thomsen \cite{Tho:HomHetOneSol}.  It follows from Theorem \ref{mainClassThm} that the homoclinic $\mathrm{C}^*$-algebra has the inductive limit structure of an ASH $\mathrm{C}^*$-algebra \cite[Lemma 5.1]{LinNiu:RationallyTAI}.  However, it would be interesting to construct an inductive limit structure for the homoclinic $\mathrm{C}^*$-algebra of a general Smale space directly from the dynamics.

\emph{Acknowledgments.} The first seeds for the ideas in this paper were planted during a visit by the first author to the Fields Institute's Thematic Programme on  Abstract Harmonic Analysis, Banach and Operator Algebras during the spring of 2014, where the second author was a postdoctoral fellow. Further collaboration was made possible by the Banach Center's sponsorship of a small workgroup in $\mathrm{C}^*$-algebras and Dynamics in March 2015. The bulk of this work was carried out in July 2015 at the University of M\"unster, where both authors were hosted during the focussed programme on $\mathrm{C}^*$-algebras. The authors thank Ian Putnam, Stuart White and Rufus Willett for useful discussions. Finally, the authors thank the participants of the conference Noncommutative Dimension Theories, held at the University of Hawaii, whose feedback on a talk given by the first listed author resulted in Theorem \ref{stableAndUnstableWithProj} and Question \ref{projQues}.

\providecommand{\bysame}{\leavevmode\hbox to3em{\hrulefill}\thinspace}
\providecommand{\MR}{\relax\ifhmode\unskip\space\fi MR }
\providecommand{\MRhref}[2]{%
  \href{http://www.ams.org/mathscinet-getitem?mr=#1}{#2}
}
\providecommand{\href}[2]{#2}

\end{document}